\documentclass[leqno,12pt]{article}
\usepackage{latexsym}
\newtheorem{thm}{Theorem}

\newtheorem{prop}[thm]{Proposition}

\newtheorem{con}[thm]{Conjecture}

\newcommand{\pf}{\noindent{\bf Proof\ \ }}
\newcommand{\finpf}{\hfill{$\Box$}\linespace}
\newcommand{\linespace}{\vspace{\baselineskip} \noindent}
\newcommand{\R}{{\bf R}}
\newcommand{\bH}{{\bf H}}
\newcommand{\Sd}{{\bf S}^d}
\newcommand{\ww}{\widehat{w}}
\newcommand{\wG}{\widehat{G}}

\newenvironment{myequation}{\setcounter{equation}{\value{thm}}
   \begin{equation}}{\addtocounter{thm}{1}\end{equation}}

\newcommand{\bmye}{\begin{myequation}}
\newcommand{\emye}{\end{myequation}}

\begin{document}
\title{\textbf{\normalsize
THE LAX CONJECTURE IS TRUE}}
\author{\textbf{A.S. Lewis\thanks{Department 
of Mathematics, Simon Fraser University, Burnaby, BC V5A 1S6, Canada. 
\texttt{aslewis\char64 sfu.ca, 
http://www.cecm.sfu.ca/\~~\hspace{-4pt}aslewis}. 
Research supported by NSERC.}\,,
P.A. Parrilo\thanks{Automatic Control Laboratory, Swiss Federal Institute
of Technology, CH-8092 Z\"urich, Switzerland.
\texttt{parrilo\char64 control.ee.ethz.ch}.}\,, and 
M.V. Ramana\thanks{Corporate Research and Development, United Airlines Inc., 
Elk Grove Village, Illinois, USA.
\texttt{motakuri\_ramana\char64 yahoo.com}. }
}}
\date{}

\maketitle


\begin{abstract}
In 1958 Lax conjectured that hyperbolic polynomials in three variables
are determinants of linear combinations of three symmetric matrices.
This conjecture is equivalent to a recent observation
of Helton and Vinnikov.
\end{abstract}

Consider a polynomial $p$ on $\R^n$ of degree $d$ (the maximum of the
degrees of the monomials in the expansion of $p$).  We call $p$
{\em homogeneous} if $p(tw) = t^d p(w)$ for all real $t$ and vectors
$w \in \R^n$:  equivalently, every monomial in the expansion of $p$
has degree $d$.  We denote the set of such polynomials by $\bH^n(d)$.  By identifying a polynomial with its vector of coefficients, we can consider $\bH^n(d)$ as a normed vector space of dimension ${n+d-1 \choose d}$.

A polynomial $p \in \bH^n(d)$ is {\em hyperbolic} with respect
to a vector $e \in \R^n$ if $p(e) \ne 0$ and, 
for all vectors $w \in \R^n$, the univariate
polynomial $t \mapsto p(w-te)$ has all real roots.  The corresponding
{\em hyperbolicity cone} is the open convex cone (see \cite{Gar59})
\[
\{ w \in \R^n : p(w-te) = 0~\Rightarrow~t > 0 \}.
\] 
For example, the polynomial $w_1 w_2 \cdots w_n$ is hyperbolic with respect
to the vector $(1,1,\ldots,1)$, since the polynomial 
$t \mapsto (w_1 - t)(w_2 - t) \cdots (w_n - t)$ has roots 
$w_1,w_2,\ldots,w_n$; hence the corresponding
hyperbolicity cone is the open positive orthant.

Hyperbolic polynomials and their hyperbolicity cones
originally appeared in the partial differential
equations literature \cite{Gar51}.  They have attracted attention
more recently as fundamental objects in modern convex optimization
\cite{Gul97,Bau01}.  Three primary reasons drive this interest:
\begin{enumerate}
\item[(i)]
the definition of ``hyperbolic polynomial'' is strikingly simple;
\item[(ii)]
the class of hyperbolic polynomials, although not well-understood,
is known to be rich --- specifically, its interior in $\bH^n(d)$ is nonempty;
\item[(iii)]
optimization problems posed over hyperbolicity cones, with linear
objective and constraint functions, are amenable to efficient
interior point algorithms.
\end{enumerate}
For more details on these reasons, see \cite{Gul97,Bau01}.

In light of the interest of hyperbolic polynomials to optimization
theorists, it is therefore natural to ask:  how general is the class
of hyperbolicity cones?  In particular, do hyperbolicity cones 
provide a more general model for convex optimization than
``semidefinite programming'' (the study of optimization problems
with linear objectives and constraints and semidefinite matrix
variables \cite{Nes94})?

We begin with some easy observations.  A rich source of examples of hyperbolicity cones are {\em semidefinite slices}, by which we mean sets 
of the form
\bmye \label{slice}
\Big\{
w : \sum_{j=1}^n w_j G_j \in \Sd_{++}
\Big\},
\emye
for matrices $G_1,G_2,\ldots,G_n$ in the space $\Sd$ of
all $d$-by-$d$ real symmetric matrices, where $\Sd_{++}$ denotes the positive definite cone.  Such cones are, in particular, 
``semidefinite representable'' in the sense of \cite{Nes94}.

\begin{prop}
Any nonempty semidefinite slice is a hyperbolicity cone.
\end{prop}

\pf
Suppose the semidefinite slice (\ref{slice}) contains the vector $\ww$.  We claim the polynomial $p$ on $\R^n$ defined by
\bmye \label{poly}
p(w) = \det \sum_j w_j G_j
\emye
is hyperbolic with respect to $\ww$, with corresponding hyperbolicity described by (\ref{slice}).  Clearly $p$ is homogeneous of degree $d$, and $p(\ww) > 0$.

Define a matrix $\wG = \sum_j \ww_j G_j \in \Sd_{++}$, and notice, for any vector $w \in \R^n$ and scalar $t$, we have
\begin{eqnarray*}
\lefteqn{
p(w-t\ww) = \det \sum_j (w_j - t\ww_j)G_j 
= \det \Big( \sum_j w_j G_j - t\wG \Big) }  \\
& & \mbox{\hspace{2cm}} = (\det \wG)
\det \Big( 
\wG^{-1/2} \Big[ \sum_j w_j G_j \Big] \wG^{-1/2} - tI 
\Big),
\end{eqnarray*}
where $I$ denotes the identity matrix.  Consequently, the univariate polynomial $t \mapsto p(w-t\ww)$ has all real roots, namely the eigenvalues of the symmetric matrix $H = \wG^{-1/2} [ \sum_j w_j G_j ] \wG^{-1/2}$, so $p$ is hyperbolic with respect to $\ww$.  Furthermore, by definition, $w$ lies in the corresponding hyperbolicity cone exactly when these roots (or equivalently, eigenvalues) are all strictly positive.  But this property is equivalent to $H$ being positive definite, which holds if and only if 
$\sum_j w_j G_j$ is positive definite, as required.
\finpf

The class of semidefinite slices is quite broad.  For example,
any {\em homogeneous cone} (an open convex pointed cone whose automorphism group acts transitively) is a semidefinite slice \cite{Chu03} (see also \cite{Fay02}).  In particular, therefore, any homogeneous cone is a hyperbolicity cone, a result first observed in \cite{Gul97}.

What about the converse?  When is a hyperbolicity cone a semidefinite slice?  How general is the class of hyperbolic polynomials of the form (\ref{poly})?

In considering a general 
hyperbolic polynomial $p$ on $\R^n$ with respect to a vector $e$,
we can suppose, after a change of variables, that 
$e = (1,0,0,\ldots,0)$ and $p(e)=1$.  Consider the first nontrivial
case, that of $n=2$.  By
assumption, the polynomial $t \mapsto p(-t,1)$ has all real
roots, which we denote $g_1,g_2,\ldots,g_d$, so for some
nonzero real $k$ we have the identity
\[
p(-t,1) = k \prod_{j=1}^d (g_j - t).
\]
By homogeneity, for any vector $(x,y) \in \R^2$ with $y \ne 0$,
we deduce
\[
p(x,y) = y^d p\Big( \frac{x}{y} , 1 \Big)
 = y^d k \prod_{j=1}^d \Big( g_j + \frac{x}{y} \Big)
 = k \prod_{j=1}^d (g_j y + x).
\]
By continuity and the fact that $p(1,0) = 1$, we see
\[
p(x,y) = \prod_{j=1}^d (g_j y + x) = \det(xI + yG)
\]
for all $(x,y) \in \R^2$, where $G$ is the diagonal matrix
with diagonal entries $g_1,g_2,\ldots,g_d$.  Thus
any such hyperbolic polynomial $p$ does indeed have
the form (\ref{poly}).

What about hyperbolic polynomials in more than two variables?
The following conjecture \cite{Lax58} proposes that
all hyperbolic polynomials in three variables are likewise
easily described in terms of determinants of symmetric matrices.

\begin{con}[Lax, 1958]
A polynomial $p$ on $\R^3$ is hyperbolic of degree $d$
with respect to the vector $e=(1,0,0)$ and satisfies
$p(e)=1$ if and only if there exist matrices 
$B,C \in {\bf S}^d$ such that $p$ is given by
\bmye \label{detp}
p(x,y,z) = \det(xI+yB+zC).
\emye
\end{con}

\noindent
An obvious consequence of this conjecture would be
that, in $\R^3$, hyperbolicity cones and semidefinite slices comprise identical classes.

A polynomial on $\R^2$ is a {\em real zero polynomial} \cite{Hel02}
if, for all vectors $(y,z) \in \R^2$, the univariate polynomial 
$t \mapsto q(ty,tz)$ has all real roots.  Such polynomials are
closely related to hyperbolic polynomials via the following
elementary result.

\begin{prop} \label{equiv}
If $p$ is a hyperbolic polynomial of degree $d$
on $\R^3$ with respect to the vector
$e=(1,0,0)$, and $p(e)=1$, then the polynomial on $\R^2$
defined by $q(y,z) = p(1,y,z)$ is a real zero polynomial of 
degree no more than $d$, and satisfying $q(0,0)=1$.  

Conversely, if $q$ 
is a real zero polynomial of degree $d$ on $\R^2$ satisfying $q(0,0)=1$,
then the polynomial on $\R^3$ defined by
\bmye \label{hom}
p(x,y,z) = x^d q\Big( \frac{y}{x},\frac{z}{x} \Big)~~~(x \ne 0)
\emye
(extended to $\R^3$ by continuity)
is a hyperbolic polynomial of degree $d$
on $\R^3$ with respect to $e$, and $p(e)=1$.
\end{prop}

\pf
To prove the first statement, note that
for any point $(y,z) \in \R^2$ and complex $\mu$,
if $q(\mu(y,z)) = 0$ then $\mu \ne 0$ and
$0 = p(1,\mu y,\mu z) = \mu^d p(\mu^{-1},y,z)$,
using the homogeneity of $p$.  So, by the hyperbolic property,
$-\mu^{-1}$ is real, and hence so is $\mu$.
The remaining claims are clear.

For the converse direction, since $q$ has degree $d$,
clearly $p$ is well-defined and homogeneous of degree $d$
and satisfies $p(e)=1$.  If $p(\mu,y,z) = 0$, then either $\mu=0$ or
$q(\mu^{-1}(y,z)) = 0$, in which case $\mu^{-1}$ and hence also $\mu$ must
be real. 
\finpf

\noindent
(Notice, in the first claim of the proposition, that the polynomial
$q$ may have degree strictly less than $d$:  consider, for example,
the case $p(x,y,z) = x^d$.)

Helton and Vinnikov \cite[p.~10]{Hel02} observe the following result, based
heavily on \cite{Vin93}.

\begin{thm} \label{vin}
A polynomial $q$ on $\R^2$ is a real zero polynomial of degree $d$ and 
satisfies $q(0,0)=1$ if and only if there exist matrices 
$B,C \in {\bf S}^d$ such that $q$ is given by
\bmye \label{det}
q(y,z) = \det(I+yB+zC).
\emye
\end{thm}

\noindent
(Notice, as in the Lax conjecture, the ``if'' direction is immediate.)

We claim that Theorem \ref{vin} is equivalent to the Lax conjecture.  
To see this, suppose
$p$ is a hyperbolic polynomial of degree $d$
on $\R^3$ with respect to the vector
$e=(1,0,0)$, and $p(e)=1$.  Then by Proposition \ref{equiv}, 
the polynomial on $\R^2$
defined by $q(y,z) = p(1,y,z)$ is a real zero polynomial of 
degree $d'\le d$, and satisfying $q(0,0)=1$. Hence by
Theorem \ref{vin}, equation (\ref{det}) holds:
we can assume $d'=d$ by replacing $B,C \in {\bf S}^{d'}$ with 
block diagonal matrices $\mbox{Diag}(B,0),\mbox{Diag}(C,0) \in \Sd$.
Then, by homogeneity, for $x \ne 0$,
\begin{eqnarray*}
\lefteqn{
p(x,y,z) = x^d p \Big(1,\frac{y}{x},\frac{z}{x}\Big)
 = x^d q \Big(\frac{y}{x},\frac{z}{x}\Big)
}  \\
 & & \hspace{2cm}
 = x^d \det \Big(I + \frac{y}{x}B + \frac{z}{x}C\Big)
 = \det(xI + yB + zC).
\end{eqnarray*}
as required.  The converse direction in the Lax conjecture is immediate.

Conversely, let us assume the Lax conjecture, and suppose
$q$ is a real zero polynomial of degree $d$ on $\R^2$ satisfying 
$q(0,0)=1$.  (The converse direction in Theorem \ref{vin}
is immediate.)  Then by Proposition \ref{equiv} the polynomial
$p$ defined by equation (\ref{hom}) is a hyperbolic polynomial 
of degree $d$ on $\R^3$ with respect to $e$, and $p(e)=1$.
According to the Lax conjecture, equation (\ref{detp}) holds, so
\[
q(y,z) = p(1,y,z) = \det(I+yB+zC),
\]
as required.
\finpf

The exact analogue of the Lax conjecture fails in general for
polynomials in $n>3$ variables.  To see this, note that the set of
polynomials on $\R^n$ of the form $w \mapsto \det \sum_j w_j G_j$
(where $G_1,G_2,\ldots,G_n \in \Sd$) has dimension at most
$n \cdot {d+1 \choose 2}$, being an algebraic image of a vector space
of this dimension.
If the degree $d$ is large, this dimension is
certainly smaller than the dimension of the set of hyperbolic
polynomials:  as we observed above, this latter set has nonempty interior
in the space ${\bf H}^n(d)$ (by a result of Nuij \cite[Thm 2.1]{Gul97}),
and so has dimension ${n+d-1 \choose d}$.

More concretely, consider the polynomial defined by $p(w) = w_1^2 -
\sum_2^n w_j^2$ for $w \in \R^n$.  This polynomial is hyperbolic of
degree $d=2$ with respect to the vector $(1,0,0,\ldots,0)$, and yet
cannot be written in the form $\det \sum_j w_j G_j$ for matrices
$G_1,G_2,\ldots,G_n \in {\bf S}^2$ if $n>3$.  To see this, choose any
nonzero vector $w$ satisfying $w_1 = 0$, and such that the first row
of the matrix $\sum_j w_j G_j$ is zero.

The question of whether all hyperbolicity cones are semidefinite slices, or, more generally, are semidefinite representable, appears open.

\medskip

\noindent
{\bf Acknowledgement}  We are very grateful to the Institute for Mathematics
and its Applications at the University of Minnesota for their hospitality
during our work on this topic.


\end{document}